\documentclass[11pt]{article}
\usepackage{amssymb}
\title{Polynomial Relations Among Characters coming from Quantum
Affine Algebras}
\author{Michael Kleber}
\date{} 
\newtheorem{thm}{Theorem}
\newtheorem{cor}[thm]{Corollary}

\newcommand{\ba}{\leftarrow}

\newcommand{\g}{{\mathfrak g}}
\newcommand{\ghat}{{\hat{\mathfrak g}}}
\newcommand{\gl}{{\mathfrak{gl}}}
\newcommand{\into}{\hookrightarrow}
\renewcommand{\l}{\ell}
\newcommand{\rank}{{\mathop{\mathrm{rank}}\nolimits}}
\newcommand{\tensor}{\otimes}
\newcommand{\w}{\omega}

\begin{document}
\maketitle
\begin{abstract}
The Jacobi-Trudi formula implies some interesting quadratic identities
for characters of representations of $\gl_n$.  Earlier work of
Kirillov and Reshetikhin proposed a generalization of these identities
to the other classical Lie algebras, and conjectured that the
characters of certain finite-dimensional representations of
$U_q(\ghat)$ satisfy it.  Here we use a positivity argument to show
that the generalized identities have only one solution.
\end{abstract}

\section{Introduction}
\label{sec_intro}

\subsection{Motivation}

The Jacobi-Trudi (or Giambelli) formula tells us that the Schur
function of an arbitrary partition can be realized as the determinant
of a matrix whose entries are homogeneous (or elementary) symmetric
functions.  In the language of representations of $\gl_n$ indexed by
Young diagrams, this says that the character of an arbitrary
representation is a determinant of a matrix whose entries are
characters for Young diagrams with a single row (or column).

Now look at representations corresponding to rectangular Young
diagrams.  The matrix coming from an $\l\times (m+1)$ rectangle
contains as minors the matrices corresponding to rectangles of sizes
$(\l-1)\times m$, $(\l+1)\times m$, $\l\times (m-1)$, and $\l\times m$
in two different ways.  The three-term Pl\"ucker relation then yields
the following identity:
\begin{equation}
\label{polyrel-gln}
Q_m(\l)^2 = Q_{m-1}(\l) Q_{m+1}(\l) + Q_m(\l-1) Q_m(\l+1)
\end{equation}
where $Q_m(\l)$ is the character associated to the $\l\times m$
rectangular Young diagram.  This beautiful identity is not as
well-known as it ought to be.

The representations whose Young diagrams are a single column are the
fundamental representations of $\gl_n$, and one might hope that a
similar picture could be constructed starting with the fundamental
representations of other Lie algebras.  Unfortunately, one can easily
check that determinants filled in with those fundamental characters do 
not give characters of actual representations.  Variants on the
Jacobi-Trudi identity for other groups do exist (see Appendix A.3 of
\cite{FH}), but they do not behave well with respect to taking minors, 
and so do not yield analogs of equation~(\ref{polyrel-gln}).

We can hope for a more satisfactory generalization, though: perhaps we
could start with some other representations of $\g$, not necessarily
irreducible, and build a set of representations made of their
determinants which do satisfy relations like
equation~(\ref{polyrel-gln}).

In 1987, Kirillov and Reshetikhin investigated certain representations
of a recently-defined quantum deformation of the universal enveloping
algebra of $\g$.  They conjectured that the analogs of fundamental
representations for this algebra satisfied a generalization of
equation~(\ref{polyrel-gln}).  The representations were $\g$-modules
as well, so they formed a good generalization of the complete $\gl_n$
picture.

In the present paper, we reverse this process.  Beginning with the
desire to generalize the $\gl_n$ picture to types $B$, $C$ and $D$ and
retain certain properties, we show that the Kirillov-Reshetikhin
solution is in fact the only one, regardless of its interpretation in
terms of quantum deformations.

\subsection{Background}

Let $\g$ be a complex finite-dimensional simple Lie algebra, $\ghat$ its
corresponding affine Lie algebra.  Because of the inclusion of quantum
enveloping algebras $U_q(\g)\into U_q(\ghat)$, any finite-dimensional
representation of $U_q(\ghat)$ is a direct sum of irreducible
representations of $U_q(\g)$.

Here we are particularly interested in the representations of
$U_q(\ghat)$ whose highest weights are multiples of one of the
fundamental weights $\w_1,\ldots,\w_n$ of $\g$, $n=\rank(\g)$.
Unfortunately, there is presently no character formula known for these
modules in general.  The decomposition into $U_q(\g)$-modules has been
explored in \cite{KR} and \cite{ChP}, and recently by the author in
\cite{kleber}.

Let $Q_m(\l)$ denote the character of a certain $U_q(\ghat)$-module
with highest weight $m\w_\l$, where $\l=1,\ldots,n$ and $m$ is a
nonnegative integer (see section~\ref{subsec_defs} for precise
definitions).  Based on a conjectural formula for the values of the
$Q_m(\l)$, these characters appear to satisfy certain remarkable
polynomial identities.  When $\g$ is simply-laced, the identities have
the form
\begin{equation}
\label{polyrel-simple}
Q_m(\l)^2 = Q_{m-1}(\l)\,Q_{m+1}(\l) + \prod_{\l'\sim\l} Q_m(\l')
\end{equation}
for each $\l=1,\ldots,n$ and $m\geq1$.  The product is taken over all
$\l'$ adjacent to $\l$ in the Dynkin diagram of $\g$.  When $\g=\gl_n$
this is just equation~(\ref{polyrel-gln}); the relations in full
generality are written down in Section~\ref{subsec_rels}.  Using these
relations, it is possible to write any character $Q_m(\l)$ in terms of
the characters $Q_1(\l)$ of the fundamental representations of
$U_q(\ghat)$.

The main result of this paper is that, for classical Lie algebras
$\g$, these equations have only one solution where $Q_m(\l)$ is the
character of a $U_q(\g)$-module with highest weight $m\w_\l$.  By this
condition, we mean we require that $Q_m(\l)$ is a positive integer
linear combination of irreducible $U_q(\g)$-characters whose highest
weights sit under $m\w_\l$.  We use the polynomial relations to write
some of the multiplicities with which the smaller representations
appear in $Q_m(\l)$ in terms of the multiplicities in the characters
$Q_1(\l)$.  The resulting inequalities determine all of the
multiplicities.

The author is grateful to N. Yu. Reshetikhin for suggestion of the 
problem and words of wisdom.  The research was partly supported by an 
Alfred P. Sloan Doctoral Dissertation Fellowship, and partly 
conducted while visiting the Research Institute for Mathematical 
Sciences (RIMS), Kyoto, Japan, thanks to the generosity of T. Miwa.

\section{Polynomial relations}

\subsection{Definitions}
\label{subsec_defs}

We let $\g$ be a finite-dimensional complex simple Lie algebra of rank
$n$, and $\ghat$ be its corresponding affine Lie algebra.  We will
concentrate on the classical families $A_n$, $B_n$, $C_n$ and $D_n$.
Choose simple roots $\alpha_1,\ldots,\alpha_n$ and fundamental weights
$\w_1,\ldots,\w_n$ of $\g$.

We will study certain finite-dimensional representations $W_m(\l)$ of 
$U_q(\ghat)$, where $m=0,1,2,\ldots$ and $\l=1,\ldots,n$.  Since 
$U_q(\g)$ appears as a Hopf subalgebra of $U_q(\ghat)$, we can talk 
about weights and characters of $U_q(\ghat)$ modules by restricting 
our attention to the $U_q(\g)$ action.  From this point of view, 
$W_m(\l)$ has highest weight $m\w_\l$.  The structure as a $U_q(\ghat)$
module is determined by Drinfeld polynomials $P_1(z),\ldots,P_n(z)$
instead of weights; the polynomials for $W_m(\l)$ are
\begin{eqnarray*}
P_\l(z) &=& \prod_{i=1}^{m} 
\left( z + \frac{(\alpha_i,\alpha_i)}{4}(m+1-2i) \right) \\
P_k(z) &=& 1, \mbox{ for } k\neq\l
\end{eqnarray*}
Chari and Pressley have also developed the notion of a $U_q(\ghat)$
module being a ``minimal affinization'' of an irreducible $U_q(\g)$
module; see \cite{ChP} for details.  In this language, our representation
$W_m(\l)$ is the unique minimal affinization of the irreducible
representation of $U_q(\g)$ with highest weight $m\w_\l$.

Let $Q_m(\l)$ denote the character of $W_m(\l)$ viewed as a
representation of $U_q(\g)$.  If $m=0$ then $W_m(\l)$ is the trivial
representation and $Q_m(\l)=1$.  The objects $W_1(\l)$ and $Q_1(\l)$ are
called the {\em fundamental} representations and characters.

Finally, let $V(\lambda)$ denote the character of the irreducible
representation of $U_q(\g)$ with highest weight $\lambda$.  We will write
characters $Q_m(\l)$ as sums $\sum m_\lambda V(\lambda)$.  Determining
the integers $m_\lambda$ is of interest in part because they are closely
related to solutions of certain Bethe equations; this is the subject of
\cite{KR} and \cite{kleber}.  We will refer to the coefficients
$m_\lambda$ as the multiplicity of $V(\lambda)$ in the sum.

\subsection{Relations}
\label{subsec_rels}

The characters $Q_m(\l)$ when $\g$ is of type $A_n$ satisfy
equation~(\ref{polyrel-gln}), known to mathematical physicists as the
``discrete Hirota relations.''  A conjectured generalization of these
relations appears in \cite{KR} for the classical Lie algebras, and
appear as the ``$Q$-system'' in \cite{Ku} for the exceptional cases as
well.  While we are only interested in the classical cases, we will
give the relations in full generality.

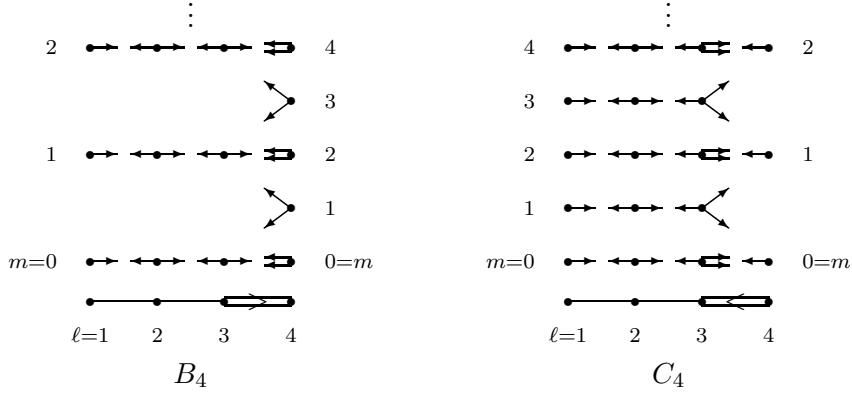
\begin{figure} \centering \setlength{\unitlength}{.035in}
\newcommand{\la}{{\vector(-1,0){4}}} \newcommand{\ra}{{\vector(1,0){4}}}
\newcommand{\ci}{{\circle*{1}}} \newcommand{\mpt}{\multiput}
\newsavebox{\lrbx}\sbox{\lrbx} {\begin{picture}(0,0)
  \mpt(0,0)(10,0){3}\ci\mpt(0,0)(10,0){2}\ra\mpt(10,0)(10,0){2}\la\end{picture}}
\newsavebox{\dd}\sbox{\dd} {\begin{picture}(0,0)
  \mpt(0,0)(10,0){4}\ci\put(0,0){\line(1,0){20}}
  \put(20,.6){\line(1,0){10}}\put(20,-.6){\line(1,0){10}}
\end{picture}}
\begin{picture}(70,60)
\mpt(20,20)(0,16){3}{\usebox{\lrbx}}
\mpt(40,20)(0,16){3}\ra
\mpt(50,20)(0,8){5}\ci
\mpt(50,20.6)(0,16){3}\la
\mpt(50,19.4)(0,16){3}\la
\mpt(50,28)(0,16){2}{\vector(-4,3){4}}
\mpt(50,28)(0,16){2}{\vector(-4,-3){4}} 
\put(20,14){\usebox{\dd}}
\put(45,14){\makebox(0,0){$>$}}
\put(20,9){\makebox(0,0){$\scriptstyle\l=1$}}
\put(30,9){\makebox(0,0){$\scriptstyle2$}}
\put(40,9){\makebox(0,0){$\scriptstyle3$}}
\put(50,9){\makebox(0,0){$\scriptstyle4$}}
\put(55,20){\makebox(0,0)[l]{$\scriptstyle 0=m$}}
\put(55,28){\makebox(0,0)[l]{$\scriptstyle 1$}}
\put(55,36){\makebox(0,0)[l]{$\scriptstyle 2$}}
\put(55,44){\makebox(0,0)[l]{$\scriptstyle 3$}}
\put(55,52){\makebox(0,0)[l]{$\scriptstyle 4$}}
\put(15,20){\makebox(0,0)[r]{$\scriptstyle m=0$}}
\put(15,36){\makebox(0,0)[r]{$\scriptstyle 1$}}
\put(15,52){\makebox(0,0)[r]{$\scriptstyle 2$}}
\put(35,3){\makebox(0,0){$B_4$}}
\put(35,58){\makebox(0,0){\vdots}}
\end{picture}
\begin{picture}(70,60)
\mpt(20,20)(0,8){5}{\usebox{\lrbx}}
\mpt(50,20)(0,16){3}\ci
\mpt(50,20)(0,16){3}\la
\mpt(40,20.6)(0,16){3}\ra
\mpt(40,19.4)(0,16){3}\ra
\mpt(40,28)(0,16){2}{\vector(4,3){4}}
\mpt(40,28)(0,16){2}{\vector(4,-3){4}} 
\put(20,14){\usebox{\dd}}
\put(45,14){\makebox(0,0){$<$}}
\put(20,9){\makebox(0,0){$\scriptstyle\l=1$}}
\put(30,9){\makebox(0,0){$\scriptstyle2$}}
\put(40,9){\makebox(0,0){$\scriptstyle3$}}
\put(50,9){\makebox(0,0){$\scriptstyle4$}}
\put(15,20){\makebox(0,0)[r]{$\scriptstyle m=0$}}
\put(15,28){\makebox(0,0)[r]{$\scriptstyle 1$}}
\put(15,36){\makebox(0,0)[r]{$\scriptstyle 2$}}
\put(15,44){\makebox(0,0)[r]{$\scriptstyle 3$}}
\put(15,52){\makebox(0,0)[r]{$\scriptstyle 4$}}
\put(55,20){\makebox(0,0)[l]{$\scriptstyle 0=m$}}
\put(55,36){\makebox(0,0)[l]{$\scriptstyle 1$}}
\put(55,52){\makebox(0,0)[l]{$\scriptstyle 2$}}
\put(35,3){\makebox(0,0){$C_4$}}
\put(35,58){\makebox(0,0){\vdots}}
\end{picture}
\caption{Product from the definition of $\cal Q$ for $B_4$ and $C_4$}
\label{fig_product}
\end{figure}

For every positive integer $m$ and for $\l=1,\ldots,n$,
\begin{equation}
\label{polyrel}
Q_m(\l)^2 = Q_{m+1}(\l)\,Q_{m-1}(\l) + \prod_{\l'\sim\l} {\cal Q}(m,\l,\l')
\end{equation}
The product is over all $\l'$ adjacent to $\l$ in the Dynkin diagram of
$\g$, and the contribution ${\cal Q}(m,\l,\l')$ from $\l'$ is determined
by the relative lengths of the roots $\alpha_\l$ and $\alpha_{\l'}$, as
follows:
\begin{equation}
\label{Qcurly}
\newcommand{\lal}{{(\alpha_\l,\alpha_\l)}}
\newcommand{\lalp}{{(\alpha_{\l'},\alpha_{\l'})}}
\newlength{\kw} \settowidth{\kw}{$k$} \newcommand{\nok}{{\hspace{\kw}}}
{\cal Q}(m,\l,\l') =
 \left\{ \begin{array}{ll}
\displaystyle Q_m(\l') &
  \mbox{if \,} \nok\lal = \nok\lalp \\ [5pt]
\displaystyle Q_{km}(\l') &
  \mbox{if \,} \nok\lal = k\lalp \\ [2pt]
\displaystyle \prod_{i=0}^{k-1} Q_{\lfloor \frac{m+i}{k} \rfloor}(\l') &
  \mbox{if \,} k\lal = \nok\lalp
\end{array} \right.
\end{equation}
where $\lfloor x\rfloor$ is the greatest integer not exceeding $x$.  We
note that in the classical cases, the product differs from the simplified
version in equation~(\ref{polyrel-simple}) only when:
$$
\begin{array}{rll}
\g=\mathfrak{so}(2n+1),& \l=n-1: & Q_m(n-2)\,Q_{2m}(n) \\
&\l=n: & Q_{\lfloor\frac{m}{2}\rfloor}(n-1)\,
        Q_{\lfloor\frac{m+1}{2}\rfloor}(n-1) \\
\g=\mathfrak{sp}(2n),& \l=n-1: &  Q_m(n-2)\,
        Q_{\lfloor\frac{m}{2}\rfloor}(n)\,
        Q_{\lfloor\frac{m+1}{2}\rfloor}(n) \\
&\l=n: & Q_{2m}(n-1)
\end{array}
$$
The structure of the product is easily represented graphically, with a 
vertex for each character $Q_m(\l)$ and an arrow from $Q_m(\l)$ 
pointing at each term of $\prod {\cal Q}(m,\l,\l')$; see 
Figure~\ref{fig_product} for $\g$ of type $B_4$ and $C_4$.  The 
corresponding picture for $G_2$ is similarly pleasing.

Finally, we can solve equation~(\ref{polyrel}) to get a recurrence relation:
\begin{equation}
\label{recurrence}
Q_{m}(\l) = 
\frac{ Q_{m-1}(\l)^2 - \prod {\cal Q}(m-1,\l,\l') } {Q_{m-2}(\l)}
\end{equation}

Note that the recurrence is well-founded: repeated use eventually writes
everything in terms of the fundamental characters $Q_1(\l)$.  This is
just the statement that iteration of ``move down, then follow any arrow''
in Figure~\ref{fig_product} will eventually lead you from any point to
one on the bottom row.  In fact, $Q_m(\l)$ is always a polynomial in the
fundamental characters, though from looking at the recurrence it is only
clear that it is a rational function.  A Jacobi-Trudi style formula for
writing the polynomial directly was given in \cite{Hi}.

The reason that characters of representations of quantum affine algebras
are solutions to a discrete integrable system is still a bit of a
mystery.

\section{Main Theorem}

\subsection{Statement of the Main Theorem}
\label{sec_statement}

The result of \cite{KR} was to conjecture a combinatorial formula
for all the multiplicities $Z(m,\l,\lambda)$ in the decomposition
$Q_m(\l)=\sum Z(m,\l,\lambda) V(\lambda)$.  We will refer to these
proposed characters as ``combinatorial characters'' of the
representations $W_m(\l)$, although the conjecture that they are
characters of some $U_q(\ghat)$ module is unproven.

\begin{thm}[Kirillov-Reshetikhin]
\label{thm_KR}
Let $\g$ be of type $A$, $B$, $C$ or $D$.  The combinatorial characters
of $W_m(\l)$ are the unique solution to equations~(\ref{polyrel})
and~(\ref{Qcurly}) with the initial data
$$
\begin{array}{llcll}
A_n:
& Q_1(\l) &=& V(\w_\l) & 1\leq\l\leq n \\
B_n:
& Q_1(\l) &=& V(\w_\l)+V(\w_{\l-2})+V(\w_{\l-4})+\cdots & 1\leq\l\leq n-1 \\
& Q_1(n)  &=& V(\w_n) \\
C_n:
& Q_1(\l) &=& V(\w_\l) & 1\leq\l\leq n \\
D_n:
& Q_1(\l) &=& V(\w_\l)+V(\w_{\l-2})+V(\w_{\l-4})+\cdots & 1\leq\l\leq n-2 \\
& Q_1(\l) &=& V(\w_\l) & \l=n-1,n
\end{array}
$$
The solutions $Q_m(\l)$ to this recurrence are all characters of
$U_q(\g)$, and the decomposition into irreducible representations is
described combinatorially in terms of generalized ``rigged
configurations.''
\end{thm}
The explicit combinatorial formula for this solution, as given in the
paper, is computationally intractable.  An effective algorithm for
computing this solution to the recurrence relations was given by the
author in \cite{kleber}.

The main result of this paper is that the specification of initial 
data in Theorem~\ref{thm_KR} is unnecessary.

\begin{thm}
\label{thm_main}
Let $\g$ be of type $A$, $B$, $C$ or $D$.  The combinatorial characters
of the representations $W_m(\l)$ are the only solutions to
equations~(\ref{polyrel}) and~(\ref{Qcurly}) such that $Q_m(\l)$ is a
character of a representation of $U_q(\g)$ with highest weight $m\w_\l$,
for every nonnegative integer $m$ and $1\leq\l\leq n$.
\end{thm}

We need only prove that any choice of initial data other than that in
Theorem~\ref{thm_KR} would result in some $Q_m(\l)$ which is not a
character of a representation of $U_q(\g)$.  The values $Q_m(\l)$ are
always virtual $U_q(\g)$-characters, but in all other cases, some contain
representations occurring with negative multiplicity.  As an immediate
consequence, we have:

\begin{cor}
\label{cor_main}
If the characters $Q_m(\l)$ of the representations $W_m(\l)$ obey the
recurrence relations in equations~(\ref{polyrel}) and~(\ref{Qcurly}),
then they must be given by the formula for combinatorial characters in
\cite{KR}.
\end{cor}

The technique of proof is as follows.  The possible choices of initial
data are limited by the requirement that $Q_1(\l)$ be a representation
with highest weight $\w_\l$.  That is, $Q_1(\l)$ must decompose into
irreducible $U_q(\g)$-modules as
$$
Q_1(\l) = V(\w_\l) + \sum_{\lambda\prec\w_\l} m_\lambda V(\lambda)
$$
Note that we require that $V(\w_\l)$ occur in $Q_1(\l)$ exactly once.  
Furthermore, we require that for every other component $V(\lambda)$ 
that appears, $\lambda\prec\w_\l$, {\em i.e.} that $\w_\l-\lambda$ is 
a nonzero linear combination of simple roots with nonnegative integer
coefficients.

We proceed with a case-by-case proof.  For each series, we find 
explicit multiplicities of irreducible representations occurring in 
$Q_m(\l)$ which would be negative for any choice of $Q_1(\l)$ other 
than that of Theorem~\ref{thm_KR}.  The calculations for series $B$, 
$C$ and $D$ are found in sections \ref{sec_Bn}, \ref{sec_Cn} and 
\ref{sec_Dn}, respectively.

When $\g$ is of type $A_n$, no computations are necessary, because 
every fundamental root is minuscule: there are no $\lambda\prec\w_\l$ 
to worry about, no other choices for initial data to rule out.  In 
fact, $Q_m(\l)$ is just $V(m\w_\l)$ for all $m$ and $\l$, and moreover 
every $U_q(\g)$ module is also acted upon by $U_q(\ghat)$, by means of 
the evaluation representation.

\subsection{Series $B_n$}
\label{sec_Bn}

Let $\g$ be of type $B_n$.  Let $V_i$ stand for $V(\w_i)$ for $1\leq
i\leq n-1$, and $V_{sp}$ for the character of the spin representation
with highest weight $\w_n$.  For convenience, let $\w_0=0$ and $V_0$
denote the character of the trivial representation.  Finally, we denote
by $V_n$ the character of the representation with highest weight $2\w_n$,
which behaves like the fundamental representations.

There are no dominant weights $\lambda\prec\w_n$, so $Q_1(n)=V_{sp}$.
The only weights $\lambda\prec\w_a$ are $0,\w_1,\ldots,\w_{a-1}$ for
$1\leq a\leq n-1$, so we write
\begin{equation}
Q_1(a) = V_a + \sum_{b=0}^{a-1} M_{a,b} V_b
\end{equation}
Our goal is to prove that the only possible values for the
multiplicities are
\begin{equation}
\label{BMs}
M_{a,b} = \left\{
\begin{array}{ll} 1, & a-b \mbox{ even} \\ 0, & a-b \mbox{ odd} \end{array}
\right.
\end{equation}

We will show these values are necessary inductively; the proof for
each $M_{a,b}$ will assume the result for all $M_{c,d}$ with
$\lceil\frac{c-d}2\rceil < \lceil\frac{a-b}2\rceil$ as well as those
with $\lceil\frac{c-d}2\rceil = \lceil\frac{a-b}2\rceil$ and
$c+d>a+b$.  (Here $\lceil x \rceil$ is the least integer greater than or
equal to $x$.) This amounts to working in the following order:
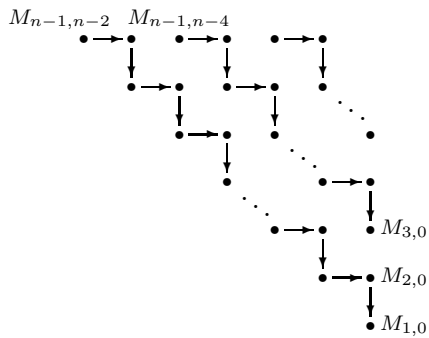
\begin{center} \setlength{\unitlength}{.25in}
\newcommand{\vx}{{\circle*{.15}}}
\newcommand{\rt}{{\vector(1,0){.6}}}
\newcommand{\dn}{{\vector(0,-1){.6}}}
\begin{picture}(8,7)
\multiput(0,6)(1,-1){7}{\vx}
\multiput(1,6)(1,-1){3}{\vx}
\multiput(5,2)(1,-1){2}{\vx}
\multiput(2,6)(1,-1){5}{\vx}
\multiput(3,6)(1,-1){2}{\vx}
\put(6,3){\vx}
\multiput(4,6)(1,-1){3}{\vx}
\put(5,6){\vx}
\multiput(0.2,6)(1,-1){3}{\rt}
\multiput(1,5.8)(1,-1){3}{\dn}
\multiput(4.2,2)(1,-1){2}{\rt}
\multiput(5,1.8)(1,-1){2}{\dn}
\multiput(2.2,6)(1,-1){2}{\rt}
\multiput(3,5.8)(1,-1){2}{\dn}
\put(5.2,3){\rt} \put(6,2.8){\dn}
\put(4.2,6){\rt} \put(5,5.8){\dn}
\multiput(3.25,2.25)(1,1){3}{\mbox{$\ddots$}}
\put(-.5,6.2){\makebox(0,0)[b]{$\scriptstyle M_{n-1,n-2}$}}
\put(2,6.2){\makebox(0,0)[b]{$\scriptstyle M_{n-1,n-4}$}}
\put(6.2,0){\makebox(0,0)[l]{$\scriptstyle M_{1,0}$}}
\put(6.2,1){\makebox(0,0)[l]{$\scriptstyle M_{2,0}$}}
\put(6.2,2){\makebox(0,0)[l]{$\scriptstyle M_{3,0}$}}
\end{picture}
\begin{minipage}[b]{2in}\raggedright
First follow the diagonal from $M_{n-1,n-2}$ to $M_{1,0}$, then the
one from $M_{n-1,n-4}$ to $M_{3,0}$, etc., ending in the top right corner
with $M_{n-1,0}$ or $M_{n-2,0}$, depending on the parity of $n$.
\vspace{12pt}
\end{minipage}
\end{center}
We show that equation~(\ref{BMs}) must hold for $M_{a,b}$, assuming it
holds for all $M_{c,d}$ which appear earlier in this ordering, by
the following calculations:
\begin{enumerate}
\item
For $M_{n-1,b}$ where $n-1-b$ is odd, the multiplicity of $V(\w_b+\w_n)$
in $Q_3(n)$ is $1-2M_{n-1,b}$,
\item
For other $M_{a,b}$ where $a-b$ is odd, the multiplicity of
$V(\w_{a+2}+\w_b)$ in $Q_2(a+1)$ is $-M_{a,b}$,
\item
For $M_{a,b}$ where $a-b$ is even:
\begin{itemize}
\item The multiplicity of $V(\w_a+\w_b)$ in $Q_2(a)$ is $2M_{a,b}-1$, and
\item The multiplicity of $V(\w_{a+2}+\w_b)$ in $Q_2(a+1)$ is $1-M_{a,b}$.
\end{itemize}
\end{enumerate}
Since all $M_{a,b}$ and all multiplicities are nonnegative integers, we
must have $M_{a,b}=0$ to satisfy the first two cases and $M_{a,b}=1$ to
satisfy the third.

The calculations to prove these claims depend on the ability to tensor
together the $U_q(\g)$-modules whose characters form $Q_m(\l)$.  A
complete algorithm for decomposing these tensors is given in terms of
crystal bases in \cite{N}.  For the current case, though, it happens that
the only tensors we need to take are of fundamental representations.
Simple explicit formulas for these decompositions had been given in
\cite{KN} before the advent of crystal base technology.

\medskip \noindent
1. $M_{n-1,b}$, $n-1-b$ odd: \\ \indent
We want to find the multiplicity of $V(\w_b+\w_n)$ in $Q_3(n)$.
Recursing through the polynomial relations, we find that
$$
Q_3(n) = Q_1(n)^3 - 2 Q_1(n) Q_1(n-1) =
Q_1(n) \left[ Q_1(n)^2-2Q_1(n-1) \right]
$$
Assuming equation~(\ref{BMs}) for $M_{n-1,b'}$ for $b'>b$ and recalling that
$Q_1(n)^2 = V_{sp}^2 = V_n+V_{n-1}+\cdots+V_0$, we need to compute
the product
$$
V_{sp}  \left[
V_n-V_{n-1}+V_{n-2}-\cdots-V_{b+1}+(1-2M_{n-1,b})V_b-\cdots \right]
$$
Since $V_{sp}V_k=\sum_{i=0}^kV(\w_i+\w_n)$, we find that the
multiplicity of $V(\w_b+\w_n)$ in the product is the desired
$1-2M_{n-1,b}$.

\medskip \noindent
2. $M_{a,b}$, $a-b$ odd, $a\leq n-2$: \\ \indent
This calculation is typical of many of the ones that will follow, and
will be written out in more detail.  We want to know the multiplicity of
$V(\w_{a+2}+\w_b)$ in $Q_2(a+1)$.  When $a\leq n-3$, we have
$$
Q_2(a+1)=Q_1(a+1)^2 - Q_1(a+2) Q_1(a)
$$
Assuming inductively that equation~(\ref{BMs}) holds for all
$M_{c,d}$ that precede $M_{a,b}$ in our ordering, we have
\begin{eqnarray*}
Q_1(a+1) &=& V_{a+1} + V_{a-1} + \cdots+V_{b}+M_{a+1,b-1}V_{b-1}+\cdots \\
Q_1(a+2) &=& V_{a+2} + V_{a} + \cdots+V_{b+1}+M_{a+2,b}V_{b}+\cdots \\
Q_1(a)   &=& V_{a} + V_{a-2} + \cdots+V_{b+1}+M_{a,b}V_{b}+\cdots
\end{eqnarray*}
To compute $Q_1(a+1)^2 - Q_1(a+2) Q_1(a)$, we note that the
$V_s V_t$ term in $Q_1(a+1)^2$ and the $V_{s+1} V_{t-1}$
term in $Q_1(a+2) Q_1(a)$ are almost identical: when $s>t$, for example,
the difference is just $\sum_{i=0}^t V(\w_i+\w_{s-t-2+i})$.  In our
case, the only $V(\w_{a+2}+\w_b)$ term that does not cancel out is the
one contributed by $M_{a,b} V_{a+2} V_b$, and the
multiplicity of $V(\w_{a+2}+\w_b)$ is $-M_{a,b}$.

When $a=n-2$ the polynomial relations instead look like
$$
Q_2(n-1)=Q_1(n-1)^2 - Q_1(n)^2 Q_1(n-2) + Q_1(n-1) Q_1(n-2)
$$
The $V_n+V_{n-2}+\cdots$ terms of $Q_1(n)^2$ behave just like the
$Q_1(a+2)$ term above.  The extra terms from
$Q_1(n-2)\left[Q_1(n-1)-V_{n-1}-V_{n-3}-\cdots\right]$ make no net
contribution, as can be seen by checking highest weights.

\medskip \noindent
3. $M_{a,b}$, $a-b$ even: \\ \indent
Calculating the multiplicity of $V(\w_{a+2}+\w_b)$ in $Q_2(a+1)$ is
similar to the above; the trick of canceling $V_s V_t$
with $V_{s+1} V_{t-1}$ works again.  The only terms remaining are
$+1$ from $V_{a+1} V_{b+1}$ and the same $-M_{a,b}$ from
$V_{a+2} M_{a,b}V_b$ as above, so the multiplicity is $1-M_{a,b}$

Likewise, calculating the multiplicity of $V(\w_a+\w_b)$ in $Q_2(a)$ we
find two contributions of $M_{a,b}$ from $M_{a,b}V_a V_b$ (in either
order) in $Q_2(a)^2$, and a contribution of $1$ from $V_{a+1} V_{b+1}$ in
$Q_1(a+1) Q_1(a-1)$, so the multiplicity is $2M_{a,b}-1$.

\subsection{Series $C_n$}
\label{sec_Cn}

Let $\g$ be of type $C_n$.  We let $V_i$ stand for $V(\w_i)$ for $1\leq
i\leq n$.  The only dominant weights $\lambda\prec\w_a$ for $1\leq a\leq
n$ are $\lambda=\w_b$ for $0\leq b<a$ and $a-b$ even, where $\w_0=0$.
(If $a-b$ is odd, then $\w_a$ and $\w_b$ lie in different translates of
the root lattice, so are incomparable.)  So we write
\begin{equation}
Q_1(a) = V_a + \sum_{i=0}^{\lfloor a/2 \rfloor} M_{a,a-2i} V_{a-2i}
\end{equation}
We will prove that in fact $M_{a,b}=0$ for all $a$ and $b$.

Again we choose a convenient order to investigate the multiplicities:
first look at $M_{a,a-2}$ for $a=n,n-1,\ldots,2$, and then all $M_{a,b}$
with $a-b=4,6,8,\ldots$.  This time the multiplicities acting as
witnesses are:
\begin{enumerate}
\item
For $M_{a,a-2}$, the multiplicity of $V(\w_{a-1}+2\w_{a-2})$ in
$Q_3(a-1)$ is $1-2M_{a,a-2}$,
\item
For $M_{a,b}$ for $a-b\geq4$, the multiplicity of $V(\w_{a-2}+\w_b)$
in $Q_2(a-1)$ is $-M_{a,b}$.
\end{enumerate}
Performing these computations requires the ability to tensor more general
representations of $\g$ than were needed in the $B_n$ case.  For this we
use the generalization of the Littlewood-Richardson rule to all classical
Lie algebras given in \cite{N}, which we summarize briefly in an Appendix.

\medskip \noindent
1. $M_{a,a-2}$: \\ \indent
We want to calculate the multiplicity of $V(\w_{a-1}+2\w_{a-2})$ in
$Q_3(a-1)$.  First we write $Q_3(a-1)$ as a sum of terms of the form
$Q_1(x)Q_1(y)Q_1(z)$, which we denote as $(x;y;z)$ for brevity.  When
$2\leq a-1\leq n-2$, we have
\begin{eqnarray*}
Q_3(a-1) &=&
(a-1;a-1;a-1)-2(a;a-1;a-2)\\
&&{}-(a+1;a-1;a-3)+(a;a;a-3)+(a+1;a-2;a-2)
\end{eqnarray*}

When $a-1$ is one of $1,2$ or $n-1$, the above decomposition still 
holds, if we set $Q_1(0)=1$ and $Q_1(-1)=Q_1(n+1)=0$.  We want to find 
the multiplicity of $V(\w_{a-1}+2\w_{a-2})$ in each of these terms.

First, $V(\w_{a-1}+2\w_{a-2})$ occurs with multiplicity 3 in the
$V_{a-1}^3$ component of $Q_1(a-1)^3$.  We calculate this
number using the crystal basis technique for tensoring representations.
Beginning with the Young diagram of $V_{a-1}$, we must choose a tableau
$1,2,\ldots,a-2,p$ from the second tensor factor, where $p$ must be be
one of $a-1$, $a$, or $\overline{a-1}$.  Then the choice of tableau from
the third tensor component must be the same but replacing $p$ with
$\overline{p}$.

Similarly, the $V_{a} V_{a-1} V_{a-2}$ component of the $(a;a-1;a-2)$
term produces $V(\w_{a-1}+2\w_{a-2})$ with multiplicity 1, corresponding
to the choice of the tableau $1,2,\ldots,a-2,\overline{a}$ from the
crystal of $V_{a-1}$.  We see that the remaining three terms cannot
contribute by looking at tableaux in the same way.

Second, $V(\w_{a-1}+2\w_{a-2})$ occurs in the $M_{a,a-2}V_{a-2} 
V_{a-1} V_{a-2}$ piece of $(a;a-1;a-2)$ and the 
$M_{a+1,a-1}V_{a-1} V_{a-2} V_{a-2}$ piece of 
$(a+1;a-2;a-2)$ as the highest weight component.  Our inductive 
hypothesis, however, assumes that $M_{a+1,a-1}=0$, and we start the 
induction with $a=n$, where the $(a+1;a-2;a-2)$ term vanishes 
entirely.

Totaling these results, we find that the net multiplicity is 
$1-2M_{a,a-2}$, and conclude that $M_{a,a-2}=0$.

\medskip \noindent
2. $M_{a,b}$ for $a-b\geq4$: \\ \indent
We want to calculate the multiplicity of $V(\w_{a-2}+\w_b)$
in $Q_2(a-1)$.  For any $2\leq a-1\leq n-1$, we have
\begin{eqnarray*}
Q_2(a-1) &=& Q_1(a-1)^2 - Q_1(a)Q_1(a-2) \\
&=& (V_{a-1}+\cdots)(V_{a-1}+\cdots) 
  - (V_a+M_{a,b}V_b+\cdots)(V_{a-2}+\cdots)
\end{eqnarray*}
where every omitted term is either already known to be 0 by induction, or
else has highest weight less than $\w_b$, so cannot contribute.  As in
the $B_n$ case, the $V_{a-1}^2$ and $V_aV_{a-2}$ terms nearly cancel one
another's contributions: their difference is just $\sum_{k=0}^{a-1}V(2\w_k)$.
Since $V(\w_{a-2}+\w_b)$ occurs in $V_{a-2}V_b$ with multiplicity 1, the
net multiplicity in $Q_2(a-1)$ is $-M_{a,b}$, and we conclude that
$M_{a,b}=0$.

\subsection{Series $D_n$}
\label{sec_Dn}

Let $\g$ be of type $D_n$.  This time we let $V_i$ stand for $V(\w_i)$
for $1\leq i\leq n-2$, and use $V_{n-1}$ for the character of the
representation with highest weight $\w_{n-1}+\w_n$.  We will not need to
explicitly use the characters of the two spin representations
individually, only their product, $V_{n-1}+V_{n-3}+\cdots$.

There are no dominant weights under $\w_{n-1}$ or $\w_n$, and so no work
to do on $Q_1(n-1)$ or $Q_1(n)$.  For $1\leq a\leq n-2$, the only
dominant weights $\lambda\prec\w_a$ are $\lambda=\w_b$ for $0\leq b<a$
and $a-b$ even; again $\w_0=0$.  (If $a-b$ is odd, then $\w_a$ and $\w_b$
lie in different translates of the root lattice, so are incomparable.)
So we write
\begin{equation}
Q_1(a) = V_a + \sum_{i=0}^{\lfloor a/2 \rfloor} M_{a,a-2i} V_{a-2i}
\end{equation}
We will show that in fact $M_{a,b}=1$ for all $a$ and $b$.

Again the proof is by induction; to show $M_{a,b}=1$ we will assume
$M_{c,d}=1$ as long as either $c-d<a-b$ or $c-d=a-b$ and $c>a$.  (This is
the same ordering used for the $B_n$ series after dropping the $M_{a,b}$
with $a-b$ odd.)  Our witnesses this time are:
\begin{itemize}
\item The multiplicity of $V(2\w_b)$ in $Q_2(a-1)$ is $1-M_{a,b}$, and
\item The multiplicity of $V(\w_a+\w_b)$ in $Q_2(a)$ is $2M_{a,b}-1$.
\end{itemize}
We must therefore conclude that $M_{a,b}=1$.  Since we only need to
tensor fundamental representations together, the explicit formulas given
in \cite{KN} are enough to carry out these calculations.

For any $\l\leq n-3$, the polynomial relations give us
$$
Q_2(\l)=Q_1(\l)^2-Q_1(\l+1)Q_1(\l-1)
$$

The multiplicity of $V(2\w_b)$ in $Q_2(a-1)$ is easily calculated
directly, since $V(2\w_b)$ appears in $V_r V_s$ if and only if
$r=s\geq b$, and then it appears with multiplicity one.  The $Q_1(a-1)$
term therefore contains $V(2\w_b)$ exactly $(a-b)/2$ times, while the
$Q_1(a)Q_1(a-2)$ term subtracts off $M_{a,b}-1+(a-b)/2$ of them.
Thus the net multiplicity is $1-M_{a,b}$.

To calculate the multiplicity of $V(\w_a+\w_b)$ in $Q_2(a)$ for $a\leq
n-3$, we once again use the trick of canceling the contribution from the
$V_s V_t$ term of $Q_1(a)^2$ with the $V_{s+1} V_{t-1}$ term of
$Q_1(a+1)Q_1(a-1)$.  The cancellation requires more attention this time,
since $V(\w_a+\w_b)$ occurs with multiplicity two in $V_s V_t$ when
$a-b\geq 2n-r-s$.  In the end, the only terms that do not cancel are the
contributions of $M_{a,b}$ from $V_a V_b$ and $V_b V_a$ in $Q_1(a)^2$ and
of $-1$ from $V_{b+1} V_{a-1}$ in $Q_1(a+1)Q_1(a-1)$.  Thus the net
multiplicity is $2M_{a,b}-1$.

Finally, if $a=n-2$ the polynomial relations change to
$$
Q_2(n-2)=Q_1(n-2)^2-Q_1(n-1)Q_1(n)Q_1(n-3)
$$
This change does not require any new work, though: $Q_1(n-1)Q_1(n)$ is
just the product of the two spin representations, which decomposes as
$V_{n-1}+V_{n-3}+\cdots$.  Since this is exactly what we wanted
$Q_1(\l+1)$ to look like in the above argument, the preceding calculation
still holds.

\appendix

\section*{Appendix: Littlewood-Richardson Rule for $C_n$}

This is a brief summary of a generalization of the Littlewood-Richardson
rule to Lie algebras of type $C_n$, as given in \cite{N}.  For our
purposes, we only need the ability to tensor an arbitrary representation
with one of the fundamental representations with highest weights
$\w_1,\ldots,\w_n$.

The representation with highest weight $\sum_{k=1}^{n} a_k\w_k$ is
represented by a Young diagram $Y$ with $a_k$ columns of height $k$.  For
a fundamental representation $V_k$, we create Young tableaux from our
column of height $k$ by filling in the boxes with $k$ distinct symbols
$i_1,\ldots,i_k$ chosen in order from the sequence
$1,2,\ldots,n,\overline{n},\ldots,\overline{2},\overline{1}$ in all
possible ways, as long as if $i_a=p$ and $i_b=\overline{p}$ then
$a+(k-b+1)\leq p$.  These tableaux label the vertices of the crystal
graph of the representation $V_k$.

Given a Young diagram $Y$, the symbols $1,2,\ldots,n$ act on it by adding
one box to the first, second,\ldots,$n$th row, and the symbols
$\overline{1},\overline{2},\ldots,\overline{n}$ act by removing one,
provided the addition or removal results in a diagram whose rows are
still nonincreasing in length.  The result of the action of the symbol
$i_a$ on $Y$ is denoted $Y\!\ba i_a$.

Then the tensor product $V\tensor V_k$, where $V$ has Young diagram $Y$,
decomposes as the sum of all representations with diagrams
$(((Y\!\ba i_1)\ba i_2)\cdots\ba i_k)$, where $i_1,\ldots,i_k$ range over
all tableaux of $V_k$ such that each of the actions still result in a
diagram whose rows are still nonincreasing in length.

\end{document}